\documentclass{amsart}
\usepackage{amssymb,amsmath,amsxtra, xcolor}
\newtheorem{theorem}{Theorem}
\newtheorem{corollary}[theorem]{Corollary}

\newtheorem{proposition}[theorem]{Proposition}

\newtheorem{example}[theorem]{Example}
\newtheorem{lemma}[theorem]{Lemma}
\newtheorem{remark}[theorem]{Remark}

\newcommand{\B}{{\mathcal B}}
\newcommand{\C}{{\mathbb C}}

\newcommand{\M}{{\mathcal M}}
\newcommand{\N}{{\mathbb{N}}}
\newcommand{\R}{{\mathbb{R}}}

\newcommand{\CC}{{\mathcal C}}
\newcommand{\F}{{\mathcal F}}
\newcommand{\HH}{{\mathcal H}}
\newcommand{\BH}{\B(\HH)}
\newcommand{\s}{\,^*\!}
\newcommand{\st}{\,^\circ}
\newcommand{\sr}{\,^*{\mathbb R}}

\newcommand{\1}{\mathbf{1}}

\newcommand{\vt}{\,\vert\,}
\newcommand{\abs}[1]{\left\vert#1\right\vert}
\newcommand{\norm}[1]{\left\Vert#1\right\Vert}
\newcommand{\re}[1]{{\rm{Re}\left(#1\right)}}
\newcommand{\w}[1]{\widehat#1}
\newcommand{\proj}[1]{{\rm{Proj}\left(#1\right)}}
\newcommand{\dis}[1]{\displaystyle#1}

\begin{document}

\baselineskip=15pt

\begin{center}
{\textbf{NONSTANDARD HULLS OF $C^{*}$-ALGEBRAS }}

\bigskip

\sc{Stefano Baratella and Siu-Ah Ng}
\end{center}

\bigskip

\begin{quote}
\small{\em } {\sc Abstract.} We study properties of $C^{*}$-algebras obtained from the \textit{nonstandard hull construction} (a generalization of the ultraproduct of $C^*$-algebras). Among others, we prove that the properties of being an infinite and a properly infinite $C^*$-algebra are both preserved and reflected by the nonstandard hull construction.  We also show that the property of being generated by mutually orthogonal projections is preserved by the nonstandard hull. Furthermore,  we indicate how a theory of weights can be developed in a nonstandard hull. Eventually we provide a few examples of standard results that follow from the nonstandard results previously obtained.
\end{quote}

\vskip 20pt

\section{Introduction} In this paper we study properties of the class of $C^{*}$-algebras that are obtained  from the \textit{nonstandard hull} construction. The nonstandard hull construction was introduced by Luxemburg in \cite{Lu}. It provides a way of obtaining a standard mathematical object (a Banach space, a Banach algebra, a $C^*$-algebra\dots) by starting from a corresponding \textit{internal} object.  Here ``internal'' is taken in the sense of A. Robinson's \textit{nonstandard analysis} and refers to an object living in a sufficiently \textit{saturated} (i.e. rich) \textit{nonstandard universe}. We refer the reader to \cite{AFHL} for the relevant definitions and for the construction of nonstandard universes.

We will outline the construction of the  nonstandard hull  of a $C^*$-algebra   in Section~\ref{sec2}. What is remarkable about the nonstandard hull construction is that it generalizes the \textit{ultraproduct} (of Banach spaces, Banach algebras, $C^*$-algebras\dots) construction, which is well known in Functional Analysis and in Operator Theory. Namely, every  ultraproduct, say,  of $C^*$-algebras can be realized as the nonstandard hull of some internal $C^*$-algebra.  Therefore we may reformulate the first paragraph by saying that we are interested in studying  properties of ultraproducts of $C^*$-algebras.

As for the paper content, in Section~\ref{sec2} we study how notions (like \textit{positivity, self adjointness} and \textit{ordering}\,) in an internal $C^*$-algebra relate to the corresponding notions in the nonstandard hull of the algebra. In this section, we also  prove a series of results showing the richness of a nonstandard hull. In particular, the nonstandard hull supports a functional calculus for Borel functions of a self-adjoint element. However, further results provide evidence that, with the exception of the finite dimensional case, the nonstandard hull of a $C^*$-algebra (or even of a von Neumann algebra) is never a von Neumann algebra.

In Section~\ref{sec3} we study the \textit{lifting} of projections from a nonstandard hull to the original internal algebra (this is important in order to attempt a classification of the nonstandard hulls of $C^*$-algebras). We prove that the properties of having an \textit{infinite} and a \textit{properly infinite} projection are both preserved and reflected by the nonstandard hull construction.  We also show that the property of being generated by mutually orthogonal projections is preserved by the nonstandard hull.

Some of the results about liftings that we present in Sections~\ref{sec2} and \ref{sec3} have been  proved   for ultraproducts of $C^*$-algebras in \cite{GH}.

In Section~\ref{sec4}, we develop a theory of \textit{weights} in a nonstandard hull. We show that any  internal \textit{$S$-continuous} weight defined on an internal $C^*$-algebra yields a weight defined on the nonstandard hull of the algebra. Moreover, we show that the latter weight always possess a certain degree of \textit{normality} (related to the degree of saturation of the underlying nonstandard universe).

In the final section, we briefly indicate how the results obtained in the previous sections can be applied to prove genuinely standard statements, i.e.  statements that contain no reference to nonstandard analysis nor to the nonstandard hull construction.

For the Functional Analysis and the Operator Theory background, we refer the reader to the classical monographs listed in the Bibliography. As for the background in nonstandard methods, \cite{AFHL} is a valuable reference.

\vskip 20pt

\section{Background}\label{sec1}

We mainly work with abstract $C^{*}$-algebras, but sometimes we restrict ourselves to concrete $C^{*}$-algebras, i.e. $C^{*}$-subalgebras of the algebra $\BH$ of bounded linear operators on  some complex Hilbert space $\HH.$ In any case, our $C^{*}$-algebras are always unitary, i.e. contains a unit $1.$

We identify $\C$ with $\{\lambda\cdot 1\vt\lambda\in\C\}.$  The space of continuous complex valued mappings defined on a topological space $X\,$ is denoted by ${\mathcal C}(X).$

Let $\M$ be a $C^{*}$-algebra. Given $X\subset\M,$ we write $\theta(X)$ for the images of $X$ under the function $\theta.$ But we write $X^{-1}$ and $\,X^n$ for the collection of invertible elements and $n$th power of elements from $X$ respectively.

We write $\sigma (a)=\{\lambda\in\C \vt (a-\lambda)\notin\M^{-1}\},$ the spectrum of $a\in\M.$
The spectral radius of $a\in\M$ is denoted by $\rho (a) =\sup \{\abs{\lambda}\vt \lambda\in\sigma (a)\,\}.$ We always have $\rho(a)\leq \norm{a}.$

In the sequel, by $*$-embedding of a $C^*$-algebra we mean an isometric $*$-homomorphism. A $*$-isomorphism is a surjective $*$-embedding.

Recall that for $a$ in a $C^{*}$-algebra,  $a$ is \emph{self-adjoint} if $a=a^*$ and \emph{normal} if $aa^*=a^*a.$ So self-adjoint elements are normal. From the classical functional calculus (for example \cite{Con}~VIII~Theorem~2.6), we have:

\begin{theorem}\label{T00}
Let $a$ be a normal element of a $C^{*}$-algebra $\M.$ Then the space ${\mathcal C}(\sigma (a))$ of continuous complex-valued functions on $\sigma (a)$ is $*$-isomorphic to a $C^{*}$-subalgebra of $\M.$ (In fact the $C^{*}$-subalgebra generated by $a$ and $1.$) \hfill $\Box$
\end{theorem}

 The set of self-adjoint elements from a $C^{*}$-algebra $\M$ is denoted by $\re{\M}.$ The set $\re{\M}$ is closed under sum, adjoint and multiplication by reals. Each $a\in\M$ has a \emph{unique} decomposition $a=a_1+ia_2,$ where $a_1, a_2\in\re{\M}.$ For normal $a$ we always have $\norm{a}=\rho (a).$

The set of positive elements is denoted by $\M_+ =\{a^2\vt a\in\re{\M}\}.$ Recall that:
\begin{itemize}
  \item $\M_+$ is a normed-closed cone with vertex at $0;$
  \item $\dis \M_+ =\{a^* a \vt a\in\M\};$
  \item  if $\M\subseteq\BH,$ then $\dis \M_+ =\{a\in\M \vt \forall \xi\in\HH\;\langle a(\xi), \xi\rangle\geq 0\, \}=\M\cap\BH_+;$
  \item for $a,b\in\M_+,$ $ab\in\M_+$ iff $ab=ba;$
  \item for any $a\in\M_+$ and $r\in \R^+,$ we have $a^r\in\M_+,$ where $a^r$ is the element given by the functional calculus (see Theorem~\ref{T00});
  \item every $\dis a\in\re{\M}$ has a unique decomposition $a=a_1-a_2,$ where $a_1, a_2\in\M_+$ and $a_1 a_2=a_2 a_1=0;\;$ together, every $a\in\M$ has a decomposition $a=(b_1-b_2)+i(c_1-c_2),\;$ where $\; b_1, b_2, c_1,c_2\in\M_+,\; c_1 c_2=c_2 c_1=b_1 b_2=b_2 b_1=0;$
  \item (polarization) for each $0\neq a\in\M$ there are unique $\abs{a}=(a^* a)^{1/2}\in\M_+$ and partial isometry $v\in\M$ such that $a=v\abs{a}.$
\end{itemize}

For $a,b\in\M$ we define $a\leq b$ iff $b-a\in\M_+.$

We let $\proj{\M}$  denote the set of projections in $\M.\,$ So $\proj{\M}\subset \M_+$ and, for a  normal element $p,\,$  $p\in\proj{\M}$ iff $\sigma(p)\subseteq\{0,1\}$ (see \cite{Bl} Corollary II.2.3.4).  Moreover,  for all $p\in\proj{\M},$
$p=0$ iff $\sigma(p)=\{0\};\,$  $p=1$ iff  $\sigma(p)=\{1\}\,$ and $p\leq 1.$

In this paper we will deal with \textit{internal} $C^{*}$-algebras, i.e. those living in some \emph{nonstandard universe}. We recall that a nonstandard universe allows to properly extend each infinite mathematical object $X$ under consideration to an object $\s X$, in a way that  $X$ and $\s X$ both satisfy  all the first order properties expressible in the formal language of set theory by means of formulas with bounded quantifiers. This  is referred to as  the \emph{Transfer Principle}.
Details of the construction of  a nonstandard universe and the methodology of A. Robinson's nonstandard analysis can be found in \cite{AFHL}.

The nonstandard extension of a standard mathematical object $\,X\,$ is denoted by $\,\s X,\,$ and this notation should not be confused with the use of $\s$ in $C^{*}$-algebras. An element from some $\,\s X\,$ is referred to as an \emph{internal} object. We identify a property with its underlying defining set, so we may speak of $\s$\,properties such as $\s$\,continuity and $\s$\,differentiability. As already said, the extension from $\,X\,$ to $\,\s X\,$ preserves the first order properties expressible in the language of set theory with bounded quantifiers  and is done simultaneous for all mathematical objects under consideration. In particular we have new objects like $\s\N\,$ (hypernatural numbers), $\sr\,$ (hyperreal numbers) and $\s\C$ (hypercomplex numbers) which still behave with respect to each other in the same formal manner as $\N,$ $\R$ and $\C.$

Elements in the set $\s\N$ are called \emph{hyperfinite}; a set counted internally by a hyperfinite number is also called hyperfinite (this is the same as $\s$finite); given $r, s\in\s\R,\,$ if $\abs{r-s}<q$ for all $q\in\R^+$ (equivalently $\abs{r-s}<1/n$ for all $n\in\N^+$), we write $r\approx s$ (and we say that $r$ and $s$ are \emph{infinitely close}). A hyperreal number  $r$ is an \emph{infinitesimal} if $r\approx 0.$ A \emph{finite} element $r$ of $\s\R$ (written $\abs{r}<\infty$) is one with $\abs{r}<n$ for some $n\in\N$; such $r$ is infinitely close to  a unique $s\in\R$ called the \emph{standard part} (in symbol: $s=\st r$). We  write $r\approx\infty$ when $r$ is positive and infinite. We use similar notions for elements in $\,\s\C\,$ or $\,\s\R^n.$ For a subset $S\subset\s\C,$ we use $\st S$ to denote $\{ \st \zeta \vt \zeta\in S\,\}.$

For some uncountable cardinal $\kappa$ sufficiently large for our purpose, we will assume throughout that the so-called \emph{$\kappa$-saturation principle} is satisfied in our construction of nonstandard objects (which is possible under a weakened form of the Axiom of Choice), namely:
\begin{quote} \emph{If $\mathcal{F}$ is a family of fewer than $\kappa$ internal sets such that $\displaystyle{\,\cap{\mathcal F}_0\neq\emptyset}$ for any finite subfamily $\mathcal{F}_0$ of $\,\mathcal{F},$ then $\displaystyle{\,\cap{\mathcal F}\neq\emptyset}.$}
\end{quote}
An internal function $\,F:\s\C\to\s\C$ is called $S$-bounded if $F(\zeta)$ is finite for all $\zeta\in\s\C.$ (Using $\omega_1$-saturation this is equivalent to having some uniform finite bound for all $F(\zeta).$)  It is called $S$-continuous if $\lambda\approx \zeta\,\Rightarrow\,F(\lambda)\approx F(\zeta)$ for all finite $\lambda, \zeta\in\s\C.$ Note the difference between $\s$continuity and $S$-continuity. If $\,F\,$ is $S$-continuous, it can be proved that there is a unique continuous $\,f:\,\C\to\C\,$ such that $\,f(\lambda)\,=\,\st\left(F(\zeta)\right)\,$ whenever $\,\zeta\approx \lambda.\,$ We write in this case $\,f\,=\,\st F$ and we say that $F$ is a \emph{lifting} of  $f.$

For readers not willing to delve into rigorous nonstandard analysis details, it suffices to regard any $\s X$ as some ultrapower $\prod_{\text{U}}X$ and internal subsets of $\s X$ as some ultraproduct $\prod_{\text{U}}X_i,$  $X_i\subseteq X,$ for some ultrafilter $\text{U}.$

\vskip 20pt

\section{Nonstandard hull of a $C^*$-algebra}\label{sec2}

In this section, $\M$ always stands for an arbitrary internal $C^*$-algebra (unless stated otherwise $\M$ is not necessarily an algebra of operators) and we construct from $\M$ a usual $C^*$-algebra, called the \emph{nonstandard hull} of $\M$ by letting
\[\w\M = \big({\rm Fin}(\M)/\!\approx\big)\,=\{\w a\,\vert\, a\in {\rm Fin}(\M)\},\] where
\begin{itemize}
\item ${\rm Fin}(\M)=\{ a\in\M\vt \norm{a}<\infty\}; $
\item for $a,b\in\M,\;$ $a\approx b$ if $\norm{a-b}\approx 0;$
\item for $a\in\M,\,$ we write $\w a=\{ x\in\M\,\vert\, x\approx a\}.\,$
\end{itemize}

We define operations on $\w\M$ by $$0=\w 0;\quad 1=\w 1;\quad \w a+\w b = \widehat{(a+b)};\quad(\w a)(\w b)=\widehat{ab};\quad (\w a)^*= \widehat{(a^*)}$$ and $\norm{\w a}=\st\norm{a}.$

\begin{lemma}\label{L200}
The operations on $\w\M$ are well defined and $\w\M$ is a $C^*$-algebra.
\end{lemma}

\begin{proof}
It is straightforward to check that sum and product on $\w\M$ are well-defined and $\w\M$ is a Banach algebra (see \cite{BN}).

As for the  map $^*$ on $\w M$, we note that if $a, b\in {\rm Fin}(\M)$ are infinitely close then $a^*\approx b^*$ and so $(\hat a)^*$ is well-defined for all $a\in {\rm Fin}(\M).$

Moreover the map $^*$ on $\w M$ is an involution and the $C^{*}$-axiom holds:
$$\norm{(\w a)^*\, (\w a)}=\norm{\widehat{a^*\, a}}\approx \norm{a^* a}=\norm{a}^2\approx \norm{\w a}^2.$$
\end{proof}

We always use $\HH$ to denote an internal complex Hilbert space.

Similar to the above, one defines the nonstandard hull $\w\HH$ of $\HH$ and proves that $\w\HH$ is a usual Hilbert space with respect to the standard part of the inner product of $\HH.$

\begin{remark} If furthermore $\M\subseteq\BH,$ then  each element $\w a\in\w\M$ can be regarded as an element of $\B(\w\HH)$ by letting $\w a(\w x)=\widehat{a(x)},$ for all $x\in\HH$ of finite norm.  (Note that $\w a(\w x)$ is well defined since $a$ is norm--finite.)
Therefore we can regard $\w\M$ as a $C^{*}$-subalgebra of $\B(\w\HH).$\hfill $\Box$
\end{remark}

\begin{remark} It can be proved that any ultrapower (in the functional analysis sense) of $C^{*}$-algebras can be realized as the nonstandard hull of an internal $C^{*}$-algebra belonging to some nonstandard universe. Hence nonstandard hulls do generalize the ultrapower construction.\hfill $\Box$
\end{remark}

\begin{proposition}\label{P03} \quad

\begin{enumerate}
\item[(i)] Let $a\in {\rm Fin}(\M).$  Then $\w a\in\w\M_+$ if and only if  there exists $c\in\M_+$ such that $\w a=\w c.$  Therefore
  $\left(\w\M\right)_+=\widehat{(\M_+)},$ and we write $\w\M_+$ unambiguously.

\item[(ii)] Let $a\in {\rm Fin}(\M).$  Then $\w a\in\re{\w\M}$ if and only if there exists $c\in\re{\M}$ such that $\w a=\w c.$ Therefore
  $\re{\w\M}=\widehat{\re{\M}}.$

\item [(iii)] Let $m_1,m_2\in \w\M.$ Then $m_1\leq m_2$ iff there are $a,b\in {\rm Fin}(\M)$ such that $m_1=\w a\;, m_2=\w b$ and $a\leq b.$ Moreover, if $m_1\in \w\M_+$ we can take $a,b\in {\rm Fin}(\M_+).$
\end{enumerate}
\end{proposition}

\begin{proof}
(i): $\w a\in \left(\w\M\right)_+$ iff $\w a=(\w b)^* (\w b)\;$  for some $\;\w b\in\w \M\;$ iff $\w a =\w c\;$ for some $\;c= b^* b\in\M_+.$

(ii): Write a self-adjoint element as the difference of two positive elements then apply (i).

(iii): One direction is trivial. Now assume $m_1\leq m_2$ in $\w\M.$  Let $m_2-m_1=\w c.$ By (i), we can take $c\in{\rm Fin}(\M_+).$ Then let $m_1=\w a$ and $b=a+c,$ so $a\leq b$ and $m_2=\w b.$ If in addition $m_1\in \w\M_+,$ we take $a\in{\rm Fin}(\M_+)$ and hence $b\in{\rm Fin}(\M_+)$ as well.
\end{proof}

\begin{lemma}\label{L202}
Let $a\in {\rm Fin}(\re{\M}).$ Then $\sigma (\w a)=\st\sigma (a).$
\end{lemma}

\begin{proof}
First we have $\st\sigma (a)\subset\sigma (\w a).\,$ To see this, let $\st\lambda\in\R.$ If $(a-\lambda)\notin\M^{-1},$ then, by \cite{BN} Theorem 5,  $(\w a-\st\lambda)=\widehat{a-\lambda}\notin \big({\w\M}\big)^{-1}.$

To prove the other inclusion, we let $\st\lambda\in\sigma (\w a )$ and show for some $\zeta \in \sigma (a)$ that $\lambda\approx \zeta.$ For non-triviality, we assume that $\lambda\notin\sigma (a).$ That is, we have $(a-\lambda)\in\M^{-1}$ but $(\w a-\st\lambda) =\widehat{a-\lambda}\notin\big(\w\M\big)^{-1}.$ Then by \cite{BN} Theorem~5, $\norm{ (a-\lambda )^{-1}}$ is infinite. By  \cite{Aup} Theorem~3.3.5  and $(a-\lambda)^{-1}\in\re{\M},$
\[{\rm dist}\big(\lambda, \sigma (a)\big) =\frac{1}{\rho \big( (a-\lambda)^{-1}\big)}=\frac{1}{\norm{ (a-\lambda)^{-1}}}.\]
Therefore ${\rm dist}\big(\lambda, \sigma (a)\big)\approx 0,$ i.e. $\lambda\approx\zeta$ for some $\zeta\in\sigma (a).$
\end{proof}

By Theorem \ref{T00}, if $a\in {\rm Fin}(\M )$ is normal, both ${\mathcal C}(\sigma (\w a))$ and  $\widehat{{\mathcal C}(\sigma (a))}$ $*$-embeds into $\w \M.$ In the case when $\w a$ is self-adjoint, the following shows that the embedding of ${\mathcal C}(\sigma (\w a))$ can be done via a lifting first to $\widehat{{\mathcal C}(\sigma (a))}.$

\begin{theorem}\label{T201} Let $a\in {\rm Fin}(\re{\M}).$ Then for every $f\in{\mathcal C}(\sigma (\w a))$ there is an $S$-continuous $F\in \s {\mathcal C}(\R )$ such that $f=\st \big( F\upharpoonright_{\sigma (a)}\big).$

In particular ${\mathcal C}(\sigma (\w a))$ $*$-embeds  into $\w \M$ via this lifting.
\end{theorem}

\begin{proof} Let $f\in {\mathcal C}(\sigma (\w a)).$ Since $\sigma (\w a)$ is compact, we apply the Tietze Extension Theorem to extend $f$ to some $\tilde{f}\in {\mathcal C}(\R )$ which has a compact support. Then let $F = \s\tilde{f}.$ Clearly $F$ is $S$-continuous. Moreover, since $\sigma (\w a)=\st\sigma (a),$ by Lemma~\ref{L202}, we have $f=\st \big( F\upharpoonright_{\sigma (a)}\big).$
\end{proof}

The closure property in the weak operator topology for a von Neumann algebra means that we can define not only continuous functions but bounded Borel functions of a normal operator. (See \cite{Aver} \S 2.6.) The following shows that bounded Borel functions of a self-adjoint operator in the nonstandard hull of an internal $C^{*}$-algebra can also be formed. This is due to the richness of $\widehat{{\mathcal C}(\sigma (a))}$ for a self-adjoint $a.$ It is a nontrivial use of saturation, because not every nonstandard hull $C^{*}$-algebra is von Neumann. (See Example \ref{nshNvna}.)

\begin{theorem}\label{T202}
Assume that $\HH$ is a $\s$separable internal Hilbert space and $\M\subseteq\BH.$ Let $m\in\re{\w\M}.$ Then for some positive regular Borel measure $\mu$ on $\sigma (m),\;$ there is a $*$-embedding of $L^\infty ( \sigma (m), \mu )$ into $\w \M$ which extends the canonical embedding of ${\mathcal C}(\sigma (m)).$
\end{theorem}

\begin{proof}
By Lemma \ref{L202} and Proposition~\ref{P03}(ii), for some $a\in\re{\M}$ we have $m=\w a$ and $\sigma(\w a)=\st\sigma(a).$

By \cite{Dou} Theorem~4.71, there exist an internal positive regular Borel measure $\nu$ on $\sigma (a)$  and an internal $*$-embedding $\gamma$ from the internal space  $ L^\infty ( \sigma (a), \nu )$ into $\BH.$ Moreover, $\gamma$ extends the embedding $\tau : {\mathcal C}(\sigma (a)) \to \M$ given by Theorem 1, where ${\mathcal C}(\sigma (a))$ denotes the internal space of $\s$continuous functions on $\sigma(a).$

We first let $\mu$ be the positive regular Borel measure on  $\st\sigma(a) = \sigma (\w a)$ given by the Loeb measure $L(\nu )$ of $\nu$ (see \cite{AFHL} Chapter 3 and the generalization for unbounded Loeb measures), i.e. $\mu = L(\nu )\circ\, {\rm{st}}^{-1},$ where ${\rm{st}}:\sigma (a)\to \st\sigma (a)$ is the usual standard mapping restricted to $\sigma (a).$

Before we define the required $*$-embedding $\pi : L^\infty ( \sigma (\w a), \mu ) \to \w\M,\,$ we consider an arbitrary $f\in L^\infty ( \sigma (\w a), \mu ) = L^\infty ( \st\sigma (a), \mu ).$ (For convenience, $f$ stands for both the function and the equivalence class it represents.)
Let $I\subset \R$ be a bounded closed interval that includes $\sigma (\w a)$ as a subset in its interior. We can regard $\nu$ and $\mu$ as measures on $\s I$ and $I$ by assigning measure $0$ to the complement of $\s I\setminus\sigma (a)$ and $I\setminus\sigma (\w a)$ respectively. Then we treat $f$ as an element in $ L^\infty ( I, \mu )$ similarly by defining $f=0$ on $I\setminus\sigma (\w a).$

Note that $\nu$ is necessarily $\s\sigma$-finite, because positive continuous functions such as the square function are in $ L^\infty ( \sigma (a), \nu )\,,$ as $a^2\in\M.$ (If $\, 0\notin \sigma (a),$ then actually $\nu$ is $\s$-finite. ) So without loss of generality, we may assume that the support of $\s f$ has $\s$finite $\nu$-measure. Now apply Lusin's Theorem (\cite{Rud}) to the $\s$compact $\s I$ and internal regular Borel measure $\nu,$ there is $\theta\in\s{\mathcal C}(I)$ so that $\norm{\theta -\s f}_{L^\infty (\s I,\nu )} \approx 0,$ where $L^\infty (\s I,\nu )$ is the internal space of $L^\infty$-functions on $\s I$ with respect to the measure $\nu.$
Therefore we have
\begin{align*}
\norm{f}_{L^\infty (\sigma(\w a),\mu )} &= \norm{f}_{L^\infty (I,\mu )}\,=\,\st \norm{\s f}_{ L^\infty (\s I,\nu )}\,=\,\st \norm{\theta}_{ L^\infty (\s I,\nu )}\\
&=\,\st \norm{\theta\upharpoonright_{\sigma (a)}}_{L^\infty (\sigma (a),\nu )}\, =\,\st \norm{\theta\upharpoonright_{\sigma (a)}}_{{\mathcal C}(\sigma (a))}\\
&=\,\st \norm{\tau(\theta\upharpoonright_{\sigma (a)})}_{\M}\,=\, \norm{\big(\tau(\theta\upharpoonright_{\sigma (a)})\big)^{\widehat{}}\,}_{\w\M }.\\
\end{align*}
In the above, the second equality is justified by the following chain of equivalences, valid for all $r\in \R,\,$ as an application of the Loeb theory:
\[\begin{array}{rl}
\norm{f}_{L^\infty (I,\mu )}\leq r \,\Leftrightarrow & \mu \{\abs{f}>r\}=0\, \Leftrightarrow\, \nu \{\abs{\s f}\geq r+n^{-1}\,\}\approx 0  \mbox{\quad for all\ } n\in\N^+ \, \Leftrightarrow\,\\  \Leftrightarrow &  \st\norm{\s f}_{L^\infty (\s I,\nu )} \leq r.
\end{array}
\]
The fifth equality follows from  $\gamma$ extending $\tau.$

Finally we define the embedding $\pi : L^\infty ( \sigma (\w a), \mu ) \to \w\M\,$ such that, for each $f\in L^\infty ( \sigma (\w a), \mu ),$ we let $\pi (f) = \big(\tau(\theta\upharpoonright_{\sigma (a)})\big)^{\widehat{}}.$

>From the above isometry and that $\gamma$ is an internal $*$-embedding, it is straightforward to check that $\pi$ is a well-defined $*$-embedding. Moreover, by Theorem~\ref{T201}, $\pi$ extends the canonical embedding of ${\mathcal C}(\sigma (m)).$
\end{proof}

In what follows we provide some simple conditions that force the nonstandard hull  a von Neumann algebra not to be a von Neumann algebra.  These conditions seem to suggest that, apart from trivial cases, the nonstandard hull of a von Neumann algebra is never von Neumann, although we do not have a proof of that.

Indeed, the nonstandard hull construction can be modified so that the property of being a von Neumann algebra is preserved.   In a further work we will investigate the so called  \textit{tracial nonstandard hull} construction (see \cite{HiOz}, or  \cite{McD} for the \textit{tracial ultraproduct}).

\begin{proposition}\label{no-vNa} Suppose $\M\subseteq\BH$ and there exists a sequence $\{ p_n\}_{n\in\mathbb N}\subset\proj{\M},$ all nonzero, with mutually orthogonal ranges (i.e. $p_n\, p_m=0$ whenever $n\neq m\,$). Then $\w\M$ is not a von Neumann algebra.
\end{proposition}

\begin{proof}
We show that $\w\M$ is not closed under the strong operator topology.
By saturation, let $\{p_n\}_{n\in\s\N}$ be an internal sequence extending $\{p_n\}_{n\in\N}.$ By Overspill there exists $M\in\s\N\setminus\N$ such that all elements of the hyperfinite sequence $\{p_N\}_{N<M}$ are nonzero projections with mutually orthogonal ranges.

Let $E_N$ be the range of $p_N.$ Then $\w p_N\in\proj{\M}$ with range $\w E_N$ and the $\w E_N$ are pairwise mutually orthogonal,  for all $N<M.$ For $n\in\N,$ let $\rho_n=\sum_{i\le n}\w p_i.$  Define $\theta\in\B (\w\HH )$ such that, for  $\w\xi\in\w\HH,\;$ $\dis \theta (\w\xi)=\lim_{n\to\infty}\rho_{n}(\w\xi).$ Note that, by $\xi\in{\rm Fin}(\HH ),$ the sequence $\{\rho_n(\w\xi)\}_{n\in\N}$ is Cauchy. Moreover $\norm{\theta}=1.$

Given finitely many $\w\xi_1, \dots, \w\xi_n\in\w\HH$ and standard $\epsilon>0,$ then for all large enough $m\in\N,$ we have $\norm{\theta (\w\xi_k)-\rho_{m}(\w\xi_k)}<\epsilon$ for all $k=1,\dots, n.$ Therefore $\theta$ belongs to the closure of $\w\M$ in $\B(\w\HH )$ under the strong operator topology.

Suppose $\theta =\w a$ for some $a\in{\rm Fin}(\M ).$ Let $N<M$ and let $\xi$ be a unit vector in $E_N.$ If $N\in\N$ then $\theta(\w\xi)=\w\xi$ and so $\norm{a(\xi)-\xi}<\frac 1{2}.$
If $N\in\s\N\setminus\N$ then $\w\xi$ is orthogonal to $\w E_n$, for all $n\in\N,$ and so it is orthogonal to the closure in the norm topology of $\oplus_{n\in\N}\w E_n.$ Hence $\theta(\w\xi)=0$ and so $\norm{a(\xi)-\xi}>\frac 1{2}.$

It follows that  the internal set $\big\{K<N: \forall x\in E_K(\norm{x}=1\rightarrow\norm{a(x)-x}<\frac 1{2})\big\}$ is equal to the external set $\N,$ which is a contradiction.

Hence $\theta$ belongs to the closure of $\w\M$ in $\B(\w\HH )$ under the strong operator topology, but $\theta\notin\w\M.$
\end{proof}

For the definition of notions such as \textit{infinite projection} used in the sequel, see \cite{Bl}~Definition III.1.3.1.

\begin{corollary} Suppose $\M\subseteq\BH$ and $\M$ contains an infinite projection. Then $\w M$ is not a von Neumann algebra.
\end{corollary}

\begin{proof} By Transfer of \cite{Bl}~ Proposition III.1.3.2, one can actually get an internal mutually orthogonal sequence $\{p_n\}_{n\in\s\N}\subset\proj{\M}\setminus\{0\}$ and Proposition~\ref{no-vNa} applies.
\end{proof}

\begin{corollary}\label{c1-novNa} Suppose $\M\subseteq\BH$ and there is $\{\rho_n\}_{n\in\mathbb N}\subset\proj{\M}$ whose sequence of ranges $\{E_n\}_{n\in\N}$ is strictly increasing. Then $\w\M$ is not a von Neumann algebra.
\end{corollary}

\begin{proof}
Define the sequence $\{p_n\}_{n\in\mathbb N}$ as follows: $p_0=\rho_0;$ $p_{n+1}=\rho_{n+1}-\rho_n.$ Then $p_{n+1}$ is the projection onto $E_{n+1}\cap E_n^\perp.$ It follows immediately that the $p_n$'s, $n>0,$ are nonzero and mutually orthogonal, hence the conclusion follows from Proposition~\ref{no-vNa}.
\end{proof}

We denote by ${\mathcal N}(T)$  the  null space of an operator $T.$  By using the \textit{functional calculus}, we can prove the following:

\begin{corollary} Let $\M$ be an internal von Neumann algebra of some $\BH.$ Suppose that there exists in $\M$ a sequence $(T_n)_{n\in\mathbb N}$ of positive operators such that ${\mathcal N}(T_{n+1})$ is strictly contained in ${\mathcal N}(T_n)$ for all $n\in\N.$ Then $\w\M$ is not a von Neumann algebra.
\end{corollary}

\begin{proof}
 Positivity of $T_n$ and the functional calculus imply that the projection $p_n$ on ${\mathcal N}(T_n)^\perp$ is in $\M,$ for all $n\in\N.$  From the assumptions it follows that ${\mathcal N}(T_n)^\perp\subset{\mathcal N}(T_{n+1})^\perp.$ Thus Corollary~\ref{c1-novNa} yields the conclusion.
\end{proof}

\begin{proposition}\label{notvna1} Suppose $\M\subseteq\BH$ and there exists some $\{\w p_n\}_{n\in\mathbb N}\subset\proj{\w\M}$ having strictly increasing ranges.  Then $\w\M$ is not a von Neumann algebra.
\end{proposition}

\begin{proof}
Define $\theta\in\B (\w\HH )$ such that, for  $\w\xi\in\w\HH,\;$
$\dis \theta (\w\xi)=\lim_{n\to\infty}\w p_n (\w\xi).$ Note that, by $\xi\in{\rm Fin}(\HH ),$ the sequence $\big\{\w p_n(\w\xi)\big\}_{n\in\N}$ is Cauchy. Moreover $\norm{\theta}=1.$

Given finitely many $\w\xi_1, \dots, \w\xi_n\in\w\HH$ and standard $\epsilon>0,$ then for all large enough $m\in\N,$ we have $\norm{\theta (\w\xi_k)-\w p_{m}(\w\xi_k)}<\epsilon$ for all $k=1,\dots, n.$ Therefore $\theta$ belongs to the closure of $\w\M$ in $\B(\w\HH )$ under the strong operator topology.

By assumption, there are $\w\eta_n\in \text{range}(\w p_n)\setminus \text{range}(\w p_{n-1}),\,\norm{\w\eta_n}=1,$ so that $\w p_k(\w\eta_m)=\w\eta_m$ if $m\leq k\in\N$ and $\w p_k(\w\eta_m)=0$  if $k< m\in\N.$

For each  $k\in\N,$  the sequence of internal sets  $\{A_{k,n}\}_{n\in\N}$ defined by $$A_{k,n}=\{\epsilon\in\sr\,\vert\, \max_{k<m<n}\norm{p_k (\eta_m )}\leq\epsilon<(n+1)^{-1}\}$$ has the finite intersection property. By saturation,
there is $\epsilon_k\approx 0$ such that
$\norm{p_k(\eta_m)}<\epsilon_k,$  for all $k<m\in\N.$

Apply saturation again, so $\big\{(p_k, \eta_k)\big\}_{k\in\N}$ extends to an infinite hyperfinite sequence $\big\{(p_k,\eta_k)\big\}_{k<M}$ with the property that, for each  $k<M,$ there exists $\epsilon_k\approx 0$  satisfying $\norm{p_k(\eta_m)}<\epsilon_k$ for all $k<m<M.$

Then  we have
\begin{equation}\theta(\w\eta_m)=\lim_{n\to\infty}\w p_n (\w\eta_m) = 0\quad \text{for all infinite}\; m<M.\tag{1}
\end{equation}
If it were $\theta =\w a$ for some $a\in{\rm Fin}(\M ),$ then
\[\norm{a(\eta_m )}\approx \norm{\theta (\w\eta_m)}=\norm{\lim_{n\to\infty}\w p_n (\w\eta_m)}=\norm{\w\eta_m}=1\quad\text{for all }\;m\in\N.\]
Then it would follow from saturation that there exists infinite $K<M$ such that
\[\norm{a(\eta_m)}\approx 1\quad\text{for all }\; m<K,\]
therefore
$\norm{\theta(\w\eta_m)}= 1$ for all $m<M,$
contradicting to (1).

Hence $\theta\notin\w\M,$ but it belongs to the closure of $\w\M$ in $\B(\w\HH )$ under the strong operator topology. That is, $\w\M$ is not von Neumann.
\end{proof}

\begin{example}\label{nshNvna}
Let $\HH =  \s\ell^2 (\N)$ be the internal Hilbert space of all internal sequences $\{a_n\}_{n\in\s\N}$ such that $\sum_{n\in\s\N}\vert a_n\vert^2\in\s\R.$  Let $e_n\in\HH$ be the usual orthonormal basis of $\HH$, i.e. $(e_n)_k =1$ if $k=n$ and $(e_n)_k =0$ otherwise. Let $p_n\in\BH$ be the projection onto the subspace generated by $e_n.$ Then, by taking $\M=\BH,$ the assumptions in Proposition~\ref{no-vNa} are satisfied and hence $\w\BH$ is not von Neumann.

In particular, the nonstandard hull of an internal $\BH$ is not necessarily a von Neumann algebra and the inclusion $\widehat{\BH }$ in $\B(\w\HH)$ could be strict.\hfill $\Box$
\end{example}

Indeed we can prove more:

\begin{proposition}\label{notvna0}
The following are  equivalent for an internal Hilbert space $\HH$:

\begin{enumerate}
\item[(i)] $\dim\HH\in\N;$
\item[(ii)] $\w\BH=\B(\w\HH);$
\item[(iii)] $\w\BH$ is a von Neumann algebra.
\end{enumerate}
\end{proposition}

\begin{proof} (i)$\Rightarrow$(ii): Let $\dim\HH=n\in\N.$ Then $\dim\w\HH=n$ and $\dim\BH=n^2.$
Hence $\dim\B(\w\HH)=n^2=\dim\w\BH.$

(ii)$\Rightarrow$(iii): Trivial.

(iii)$\Rightarrow$(i): By Proposition~\ref{notvna1}, it follows from (iii) that $\w\BH$ does not contain any infinite family of projections with strictly increasing ranges. Therefore $\w\HH$ must be finite dimensional, hence $\dim\HH\in\N.$
\end{proof}

\begin{example}\label{nshevnN}
Let $\mu$ be a $\s\sigma$-additive internal probability measure on a set $X$ such that the corresponding measure algebra is not standard finitely generated.

Therefore we have some internal disjoint $X_n\subset X,\, n\in\N,\,$ such that $\mu (X_n )>0.$

Let $\HH=L^2 (X,\mu).$ View $\M=L^\infty (X,\mu)$ as a multiplication algebra on $\HH.$ Thus, by classical results, $\M$ is an internal von Neumann subalgebra of $\BH.$

But the sequence $\{p_n\}_{n\in\N},$ where $p_n=\1_{X_n}\in\proj{\M},$ satisfies the assumptions of Proposition~\ref{no-vNa}, therefore $\w\M$ is not von Neumann.\hfill $\Box$
\end{example}

Notice that, for any $\s\sigma$-additive probability $\mu$ on an internal set $X,$ if we let $L(\mu)$ denote the Loeb measure, then by classical results, as a multiplication algebra, $L^\infty (X, L(\mu ) )$ is a von Neumann subalgebra of $\B (L^2 (X, L(\mu ) )),$ in fact, a maximal abelian von Neumann subalgebra.

On the other hand, when $L^\infty (X,\mu )$  is regarded as a subspace of $L^2 (X,\mu ),$ we have $\widehat{L^\infty (X,\mu )}=L^\infty (X, L(\mu ) ).$
(This is so, for if we let $F: X\to\s\C\,$ and $\s\norm{F}_\infty =r <\infty,$ then $\forall s> r\quad \mu\big(\{ F\geq s\}\big)\approx 0.\,$ So $L(\mu )\big(\{ F\geq s\}\big)\approx 0,\,$ thus $\norm{\w F}_\infty =\st r <\infty.$ For the other inclusion, one uses liftings of the functions.) The fact that $\widehat{L^\infty (X,\mu )}=L^\infty (X, L(\mu ) )\,$ does not cause a contradiction in the above example, since the nonstandard hull $\w\M=\widehat{L^\infty (X,\mu )}$ in the above example is taken under the norm of $\B(L^2 (X, \mu )).$

\vskip 20pt

\section{Projections in the nonstandard hulls}\label{sec3}

As before $\M$ always denotes an internal unitary $C^*$-algebra.

We begin with a series of \textit{lifting} lemmas of some relations on $\proj{\w\M}$ that will be used later to prove that some relevant properties of $C^*$-algebras are preserved and reflected by the nonstandard hull construction.

Recall that two projections  $p, q,$ are \textit{Murray-von Neumann equivalent} (notation: $p\sim q$) if for some partial isometry $u\,$ from the same C*-algebra that $u^* u = p$ and $u u^* = q.\,$ We say that $p$ is \textit{subordinate} to  $q$ (notation: $p \precsim q)$ if there is  a partial isometry $u$ from the same C*-algebra that $u^*u=p$ and $uu^*\le q.$

\begin{theorem}\label{plift01}\quad

\begin{enumerate}
  \item [(i)]  Let $\w v\in \w \M$ be a partial isometry. Then there is a partial isometry $u\in\M$ such that $u\approx v;$
  \item[(ii)] For any $\w p\in\proj{\w\M}$ there is $q\in\proj{\M}\,$ such that $\w p=\w q;$
  \item [(iii)] Given $\w {p_0},\w{q_0}\in\proj{\w\M}\,$with $\,\w {p_0}\sim\w {q_0},\,$ there are $p,q\in\proj\M$ such that
$\w p=\w{p_0},\,$ $\w q=\w{q_0}\,$ and $\,p\sim q.$
\end{enumerate}
\end{theorem}

\begin{proof}

(i):  Since $(\w v)^*\,\w v\in\proj{\w\M}\,$ and $\sigma\big((\w v)^*\,\w v \big) =\st \sigma (v^* v),\,$there is $\epsilon\approx 0$ such that
$\dis \sigma (v^* v)\subset [-\epsilon, \epsilon]\cup [1-\epsilon, 1+\epsilon].$

Let $\phi \in {\mathcal C}(\sigma (v^*v))\,$ such that $\phi (t)=0$ on $[-\epsilon, \epsilon]\cap \sigma (v^* v)$ and $\phi (t) = t^{-1/2}$ on $[1-\epsilon, 1+\epsilon]\cap \sigma (v^* v).$

Since $v^* v\in\re{\M },\,$ we have $\phi (v^* v)\in \M.\,$ Note that $\big(\phi(v^* v)\big)^* =\overline{\phi}(v^* v) =\phi (v^* v),\,$ so in fact $\phi(v^* v)\in\re{\M }.\,$

Define $u=v\,\phi(v^* v).\,$

Let $f=\text{id}_{\sigma(v^* v)}\,\phi^2\in {\mathcal C}(\sigma (v^*v)).\,$ So $f(t)=0$ on $[-\epsilon, \epsilon]\cap \sigma (v^* v)$ and $f(t) =1$ on $[1-\epsilon, 1+\epsilon]\cap \sigma (v^* v).$

Therefore $u^* u = \phi(v^* v)^* v^* v\phi(v^* v)= \phi(v^* v) v^* v\phi(v^* v)=f(v^* v).$

Moreover, $\sigma (u^* u) = f\big(\sigma (v^* v) \big)\subset\{0,1\}.\,$ It follows from \cite{Bl}~Corollary II.2.3.4 that  $u^* u\in\proj\M,\,$ i.e. $u$ is a partial isometry.

Since
$\,\dis \norm{\phi (v^* v)-v^* v}=\norm{\phi-{\rm{id}}_{\sigma (v^* v)}}_{C(\sigma(v^* v))}\leq(1-\epsilon)^{-1/2}-1+\epsilon\,\approx\,0,$
we have $\phi (v^* v)\approx v^* v.\,$

On the other hand, by $\w v\,$ being a partial isometry, we have $\w v =\w v (\w v)^*\w v\,$ (\cite{Bl}~II.2.3.5), i.e. $ v v^* v \approx v.$

Therefore $u=v\,\phi(v^* v)\,\approx v v^* v\,\approx\, v\,$ and (i) is proved.

(ii): Use a proof similar to that of (i). Alternatively, note that (ii) is a special case of (iii) which we are going to prove next.

(iii): Suppose $\w {p_0},\w{q_0}\in\proj{\w\M}\,$ and $\,\w {p_0}\sim\w {q_0}.\,$ So for some partial isometry $\w v\in\w \M,\,$ $(\w v)^*\, \w v= \w {p_0}\,$ and
$\w v\, (\w v)^* = \w {q_0}.\,$ i.e. $v^* v\approx p_0\,$ and $v v^*\approx q_0.\,$ By (i), we let $u\in\M$ be a partial isometry such that $u\approx v.$ Define $p=u^* u \in\proj\M,\,$ so $p\approx p_0.$ By \cite{Bl}~II.2.3.5, $u^*$ is also a partial isometry. Define $q=u u^*\in\proj\M, $ then $q\approx q_0.$

Therefore $p\approx{p_0},\,$ $q\approx{q_0}\,$ and $\,p\sim q.$
\end{proof}

\begin{lemma}\label{plift02} Let ${p_0},q\in\proj{\M}\,$with $\,\w {q}\leq \w {p_0}.\,$ Then there is $p\in\proj\M$ such that
$\w p=\w{p_0},\,$ and $\,q\leq p.$

Consequently, if $\w a, \w b\in\proj{\w\M}$ are such that $\w b \precsim \w a,$ then there are $p,q\in\proj{\M},\;$ $q\precsim p$ and $\w a =\w p,\; \w b =\w q.$
\end{lemma}

\begin{proof}
By \cite{Bl}~II.3.3.1, $\w p_0 \w q =\w q,\,$ i.e. $p_0 q\approx q.$

Then by \cite{Bl}~II.3.3.5, for any $n\in\N,\,$ there is $p_n\in \proj{\M}\,$ $q\leq p_n$ and $\norm {p_n -p_0}<1/n.$ Extend $\{ p_n\}$ to an internal sequence, let $p=p_N\in\proj{M},$ where $N$ is infinite, then $q\leq p$ and $p\approx p_0.$

Now if $\w a, \w b,\w c\in\proj{\w\M}$ are such that  $\w b \sim \w c \leq \w a,$ then by Theorem~\ref{plift01}\,(iii),  for some $q, q_0\in\proj{\M},$ $\w b=\w q,\,$ $\w c =\w q_0$ and $q\sim q_0.\,$ By Theorem~\ref{plift01}\,(ii), let $p_0\in\proj{\M}$ such that $\w a =\w p_0.$ So, by the above, there is some $p\in\proj{\M}$ such that $b\approx q\sim q_0\leq p\approx a.$

\end{proof}

\begin{lemma}\label{plift03} Let $p,q\in\proj{\M}.\,$
\begin{enumerate}
  \item [(i)]  If $p\sim 0,$ then $p=0.$
  \item [(ii)] If $p\approx q,$ then $p\sim q\,$ and $1-p\sim 1-q.$
  \item [(iii)] If $q\approx p\leq q,\,$ then $p=q.$
  \item [(iv)] If $p\approx 1$ then $p=1.$
  \item [(v)] If $p\approx 0$ then $p=0.$
\end{enumerate}
\end{lemma}

\begin{proof}
(i): Let $u\in\M$ with $u^* u=p,\,uu^*=0.$ Then $p=p^2=u^* u u^* u= u^* 0 u =0.$

(ii): Apply \cite{Bl}~II.3.3.4.

(iii): From \cite{Bl}~II.3.3.1, we have $0\approx q-p\in\proj{\M}.$ Hence $q-p=0$ by (i) and (ii).

(iv): Since $p\leq 1$ always holds, we have $1\approx p\leq 1.\,$ So $p=1$ by (iii).

(v): Apply (iv) to $1-p.$

\end{proof}

In the sequel, we prove some properties of projections in an internal $C^*$-algebra   $\M\subseteq\BH,$ where $\HH$ is a Hilbert space with inner product $\langle\ ,\ \rangle.$ Let
$V, W$ be closed subspaces of $\HH.$ We use the following notation:
\begin{enumerate}
\item $V\le W$ if $V$ is a subspace of $W;$
\item $V\approx W$ if for all $v\in\mbox{Fin}(V)$ there exists $w\in\mbox{Fin}(W)$ such that $v\approx w$ and conversely;
\item $V\dagger W$ if $\langle v,w\rangle\approx 0$ for all $v\in\mbox{Fin}(V), w\in\mbox{Fin}(W);$
\item $V\preceq W$ if $V\approx U,$ for some $U\le W.$
\end{enumerate}

\bigskip

\begin{lemma}\label{l1} The following are equivalent for $p,q\in\proj\M$:

\begin{enumerate}
\item $pq\approx 0;$
\item $\mbox{ran}(p)\dagger \mbox{ran}(q).$
\end{enumerate}
\end{lemma}

\begin{proof}
\item[$(1)\Rightarrow(2)$] By contraposition: let $v\in\mbox{ran}(p), w\in\mbox{ran}(q)$ be unit vectors such that $\langle v,w\rangle\not\approx 0.$ Since $\norm{pq(w)}=\norm{p(w)}\ge\vert\langle v,w\rangle\vert^{\frac 1 2},$ it follows that $\norm{pq}\not\approx 0.$
\item[$(2)\Rightarrow(1)$] Suppose there exists a unit vector $v\in\HH$ such that $pq(v)\not\approx 0.$ Write $q(v)=v_1+v_2,$ with $v_1\in\mbox{ran}(p)$ and $v_2\in\mbox{ran}(p)^\perp.$ Then $v_1\not\approx 0$ and so $\langle v_1, q(v)\rangle=\norm{v_1}^2\not\approx 0.$
\end{proof}

\begin{lemma}\label{l2} Let $V, W$ be closed subspaces of $\HH$ such that $V\approx W.$  Then $\langle u,w\rangle\approx 0$ for all $u\in\mbox{Fin}(V^\perp)$ and all $w\in\mbox{Fin}(W).$
\end{lemma}
\begin{proof} Let $u\in\mbox{Fin}(V^\perp)$ and $w\in\mbox{Fin}(W).$ Pick $v\in V$ such that $v\approx w.$ Then $$0\approx\langle u,v-w\rangle=\langle u,v\rangle-\langle u,w\rangle= -\langle u,w\rangle.$$
\end{proof}

\begin{lemma}\label{l3}
The following are equivalent for $p,q\in\proj\M$:

\begin{enumerate}
\item $\mbox{ran}(p)\approx \mbox{ran}(q);$
\item $p\approx q.$
\end{enumerate}
\end{lemma}
\begin{proof}
\item[$(1)\Rightarrow(2)$] Let $v\in\mbox{Fin}(\HH).$ Write $v=v_1+v_2,$ with $v_1\in\mbox{ran}(p)$ and $v_2\in\mbox{ran}(p)^\perp.$ Let $w\in\mbox{ran}(q)$ be such that $v_1\approx w.$ Then $q(v)\approx w+q(v_2).$ We claim that $q(v_2)\approx 0.$ For, if not then $\langle v_2, q(v_2)\rangle\not\approx 0,$ contradicting to Lemma~\ref{l2}.

Therefore $p(v)=v_1\approx w\approx q(v).$ Since the supremum of an internal set of infinitesimals is infinitesimal, we conclude that $\norm{p-q}\approx 0.$
\item[$(2)\Rightarrow(1)$] Straightforward.
\end{proof}

\begin{corollary} Assume further that $\M=\B(\HH)$ and let $p,q\in\proj\M$ be such that $pq\approx 0.$ Then there exist $p',q'\in\proj\M$ such that $p\approx p',$ $q\approx q'$ and $p'q'=0.$
\end{corollary}

\begin{proof} By Lemma~\ref{l1} and Lemma~\ref{l3}, it suffices to  prove that,  for all $V, W$ closed subspaces of $\HH$ with $V\dagger W,$ there exist closed subspaces $U, Z$ of $\HH$ such that $U\approx V;$ $Z\approx W$ and $U\perp Z.$

Let $p$ be the projection on $W.$ From $V\dagger W$ it follows that $p(v)\approx 0$  for all norm-finite $v\in V.$

The set  $\{v-p(v):\, v\in V\}$    is a linear subspace of $\HH.$  Let $U$ be its closure.  Clearly $U \perp W.$ Moreover, if $u\in U$ is a norm-finite vector then $u\approx v-p(v)$ for some $v\in V.$ Since $p(v)\approx 0,$ then $u\approx v.$ Conversely, if $v\in V$ is a norm-finite vector, then $v\approx v-p(v)\in U.$ Therefore $U\approx V.$

Finally, it  suffices to take $Z=W.$

\end{proof}

\begin{proposition}\label{p5} Let $p, q$ be projections in $\M$ with $pq\approx qp.$
Then $pq\approx r,$ where $r$ is the projection on  the closure of ${\mbox{ran}}(pq).$
\end{proposition}

\begin{proof} From $pq\approx qp$ it follows at once that $\w p\w q$ is a projection in $\w\M.$ By Theorem 17(ii), there exists a projection $r'\in\M$ such that $r'\approx pq.$
It is easy to see that ${\rm ran}(r^\prime)\approx {\rm ran}(pq)$ and ${\rm ran}(r)\approx {\rm ran}(pq), $ hence ${\rm ran}(r)\approx{\rm ran}(r^\prime)$ and Lemma~\ref{l3} applies.
\end{proof}

We now go back to the abstract setting and we deal with
\textit{infinite} and \textit{properly infinite}  $C^*$-algebras (we refer the reader to \cite{Bl}\,III.1.3.1 for the relevant definitions).

We prove that the properties of being an infinite and a properly infinite $C^*$-algebra are preserved and reflected by the nonstandard hull construction.

\begin{theorem}\label{finite}\quad

\begin{enumerate}
  \item [(i)] $\M$ is infinite iff $\w \M$ is infinite.
  \item [(ii)] $\M$ is properly infinite iff $\w \M$ is properly infinite.
\end{enumerate}
\end{theorem}

\begin{proof}
(i): Let $\M$ be infinite. So for some $p\in\proj {\M}$ we have $1\neq p \sim 1.$

Then by Lemma~\ref{plift03} (iv), $1\not\approx p,$ hence $1\neq\w p\sim 1\,$ and thus $\w \M$ is infinite.

Conversely, let $\w\M$ be infinite. Then $1\neq\w p\sim 1\,$ for some $\w p\in\proj {\w\M}.$ By Theorem~\ref{plift01}(iii) we can assume that $p \in\proj{\M}$  and that there exists $q \in\proj{\M}$ such that $p\sim q\approx 1.$ So $q=1$ by Lemma~\ref{plift03}(iv). Therefore we have $1\neq p$ and $p\sim 1$ (the latter by \cite{Bl} II.3.3.4), i.e. $\M$ is infinite.

(ii): Let $\M$ be properly infinite. So there are $p,q\in\proj{\M}$ such that $pq=0$ and $p\sim 1\sim q.$ Then clearly $\w p\,\w q=0$ and $\w p\sim 1\sim\w q,$ i.e. $\w \M$ is properly infinite.

Conversely, let $\w\M$ be properly infinite. Then for some $\w p,\w q\in\proj{\w\M},$ we have $\w p\,\w q=0$ and $\w p\sim 1\sim \w q.\,$ Moreover, we have $(1-\w p-\w q)^2 = 1-\w p-\w q,$ hence $(1-\w p-\w q)\in\proj{\w\M}.\,$ From $0\le 1-\w p-\w q$ It follows that  $\w p \leq 1-\w q\,$  and, by Theorem~\ref{plift01}(iii), Lemma~\ref{plift02} and Lemma~\ref{plift03}(iv), we can assume without loss of generality that $p,q\in\proj{\M},\,$  $1\sim p\leq 1-q\,$ and $q\sim 1.$

Then we have $0\leq (pq)^*(pq)=qpq.$ Moreover $qpq\leq 0,$ by \cite{Bl} II.3.1.8.
Therefore $(pq)^*(pq)=0$ and so $pq=0,$ by the $C^*$-axiom.
\end{proof}

In fact Theorem~\ref{finite}(i) can be generalized as follows.

\begin{theorem}\label{liftinfinite}
Let $p\in\proj{\M}.$ Then $p$ is infinite iff $\w p$ is infinite.
\end{theorem}

\begin{proof}
$(\Rightarrow$): Suppose $p$ is infinite, then for some $q\in\proj{\M},\; p\sim q\lneqq p.$ Since $p\neq q,$ Lemma~\ref{plift03}(iii) implies that $p\not\approx q\, $ i.e. $\w p\neq \w q.$

Therefore $\w p\sim\w q\lneqq \w p,$ i.e. $\w p$ is infinite.

$(\Leftarrow$): Let $\w p$ be infinite. Then there is a partial isometry $\w u\in\w \M$ so that $(\w u)^*\,\w u=\w p$ and $\w u\,(\w u)^*\lneqq \w p.\,$ By Theorem~\ref{plift01} (i), we can assume that $u\in\M$ is a partial isometry.

Then $u^*\, u\approx p$ and $u\,u^*\lnapprox q\approx p$ for some $q\in\proj{\M}$ by Lemma~\ref{plift02}.

By \cite{Bl}~II.3.3.4, there is a unitary $v\in\M$ such that $v q v^*=p.\,$ Let $w=vu.$ So $w$ is a partial isometry and
\[p\,\approx\, u^*\, u\,=\,w^*\, w\,\sim\,w\, w^*\, \lneqq\, v q v^*\,=\,p,\]
where the last inequality follows from $v$ being unitary and \cite{Bl}~II.3.1.5.

Now by Lemma~\ref{plift03} (ii), we have
\[p\,\sim\,w^*\,w\,\sim w\, w^*\, \lneqq\, p,\]
therefore $p$ is infinite.
\end{proof}

\bigskip

A standard $C^{*}$-algebra is called a \textit{$P^{*}$-algebra} if each self-adjoint element is generated by mutually orthogonal projections, i.e. it is a norm-limit of elements of the form $\dis \sum_{i=0}^N \alpha_i p_i,\,$ where $\alpha_i\in\C\,$ (in fact $\alpha_i\in\R\,$) and the $p_i$'s are projections such that $p_1 p_j=0$ whenever $i\neq j.$ We also consider internal $P^{*}$-algebras.

$P^{*}$-algebras are stable under the nonstandard hull construction:

\begin{theorem}
If $\M$ is an internal $P^{*}$-algebra, then $\w\M$ is a $P^{*}$-algebra.
\end{theorem}

\begin{proof}
Suppose $\M$ is an internal $P^{*}$-algebra and $\w a\in\w \M$ is self-adjoint. By Proposition~\ref{P03}\,(ii),  we can assume that $a\in{\rm Fin}(\re{\M }).\,$ So for some hyperfinite family of mutually orthogonal projections $p_i\in\proj\M$ and $\alpha_i\in\s\C\,$ we have $\dis a\approx \sum_{i=0}^N \alpha_i p_i.\,$

Without loss of generality, we assume that $\dis a= \sum_{i=0}^N \alpha_i p_i.\,$

Noticing that $\sigma (\alpha_i p_i )=\alpha_i \sigma (p_i )=\{0,\alpha_i\},\,$ and so, by \cite{Aup}~p.66~Ex.9,
\[\sigma(a)\,\subset\, \bigcup_{i=0}^N\sigma (\alpha_i p_i )\,=\,\{\alpha_i\,\vert\, 0\leq i\leq N\,\}\cup\{ 0\}.\]
Allow $\alpha_i=0$ if necessary, we can assume that $\sigma(a)\,=\,\{\alpha_i\,\vert\, 0\leq i\leq N\,\}.$

By Lemma~\ref{L202}, $\st\sigma (a)=\sigma (\w a ),$ so $\sigma (a)$ is $S$-bounded in $\s\R\,$ and each $\st\alpha_i\in\R.$ By re-indexing, we let $\alpha_i \le\alpha_{i+1}.$

Let $n\in\N,$ define $\dis a_n =\sum_{i=0}^N \beta_i p_i \in{\rm Fin}(\re{\M }),$ where $\beta_0=\alpha_0$ and for $i>0,$
\[\beta_i\,=\, \left\{
                 \begin{array}{ll}
                   \beta_{i-1}, & \hbox{if $\,\abs{\beta_{i-1}-\alpha_i}\leq n^{-1};$} \\
                   \alpha_i, & \hbox{otherwise.}
                 \end{array}
               \right.
\]
Since $\dis a-a_n=\sum_{i=0}^N (\alpha_i-\beta_i) p_i,\,$ by the same argument above,
\[\sigma(a-a_n) =\{\alpha_i-\beta_i\,\vert\, 0\leq i\leq N\,\}.\]
On the other hand, we have $\dis \abs{\alpha_i-\beta_i}\leq 1/n\,$ and since $a-a_n$ is self-adjoint,
\[\norm{a-a_n}\,=\,\rho(a-a_n) \leq \frac{1}{n}.\]
Consequently, $\dis \lim_{n\to\infty} \norm {\w a-\w a_n}=0.$

For a fixed $n\in\N,$ in the above $\dis a_n =\sum_{i=0}^N \beta_i p_i,$ the distinct elements from the sequence $\{\beta_i\}$ are increasing with step size $>n^{-1}.$ Therefore $\dis\{\beta_i\,\vert\, 0\leq i\leq N\,\}$ is a finite set by the $S$-boundedness of $\sigma (a)$ in $\s\R.$ Let $k$ be the cardinality of $\dis\{\beta_i\,\vert\, 0\leq i\leq N\,\}.$

Let $m(0)=0,\,m(k)=N+1,\,\gamma_0=\alpha_0\,$ and $q_0=p_0.\,$ For $0<i<k,$ let $m(i)$ be the least $j$ such that $\beta_j >\beta_{m(i-1)},\,$ $\gamma_i=\beta_{m(i)}\,$ and $\,\dis q_i =\sum_{j=m(i)}^{m(i+1)-1}p_j.$

Then $\dis a_n =\sum_{i=0}^{k-1} \gamma_i q_i.\,$ Observe that $q_i\in\proj\M,\,$ so $\w q_i\in\proj{\w\M}.\,$ Moreover the $\w q_i$ are mutually orthogonal. Hence $\w a$ is the norm-limit of the $\dis \sum_{i=0}^{k-1} \st\gamma_i \w{q_i},\,$  a finite linear combinations of mutually orthogonal
projections from $\w \M.$

Therefore $\w \M$ is a $P^{*}$-algebra.
\end{proof}

In $\M,\,$ let $E$ be the spectral measure of some $a\in\re{\M }.$ Then $\dis a=\int_{\sigma (a)} t\, dE_t$ is the norm limit of elements of the form $\dis \sum_{i=0}^{n} t_i E([t_i, t_{i+1})),\,$ where $t_i < t_{i+1}$ and $\sigma (a)\subset [t_0, t_{n+1}].\,$ Each $E([t_i, t_{i+1}))\,$ is a projection and in fact $E([t_i, t_{i+1}))\in\proj\M\,$ if $\M$ is a von Neumann algebra. The same applies for a standard $C^{*}$-algebra instead of $\M.$

So von Neumann algebras are $P^{*}$-algebras. It follows then:

\begin{corollary}
If $\M$ is an internal von Neumann algebra, then $\w \M$ is a $P^{*}$-algebra.\hfill $\Box$
\end{corollary}

On the other hand, by Proposition~\ref{no-vNa}, we have;

\begin{corollary}
Let $\M$ be an internal $P^{*}$-algebra. Then $\w\M$ is von Neumann iff $\w\M$ is finite dimensional.\hfill $\Box$
\end{corollary}

In particular, von Neumann algebras form a proper subclass of $P^{*}$-algebras. Of course $P^{*}$-algebras forms a proper subclass of $C^{*}$-algebras

\vskip 20pt

\section{Noncommutative Loeb theory}\label{sec4}

In this section we deal with $C^{*}$-algebras $\M\subseteq\BH.$ For convenience we first work with a standard $C^{*}$-algebra $\M.$

\begin{proposition}\label{P01}
If $a\in\M_+$ and $\norm{a}\leq 1,$ then $(1-a)\in\M_+.$
\end{proposition}

\begin{proof}
Let $x\in\HH$ and $a(x)=\lambda x+y,$ where $y\in\langle x\rangle^\perp.$ From the hypothesis it follows at once that $\lambda\in[0,1].$ Hence $\langle(1-a)(x),x\rangle= (1-\lambda)\langle x, x\rangle\ge 0,$ for all $x\in\HH.$
\end{proof}

A function $\theta :\M_+\to [0,\infty ]$ is a \emph{weight} if it is additive and positively homogeneous, i.e. $\theta (\lambda a) =\lambda \theta (a)$ for any $\lambda\in [0,\infty ),$ with the convention that $0\cdot\infty=0.$

\begin{proposition}\label{P02} Let $\theta :\M_+\to [0,\infty ]$  be a weight. Then   $\theta (a)\leq\theta (1)\,\norm{a}$ for all  $a\in\M_+.$ In particular, if $\theta(1)=0$ then $\theta(a)=0$ for all $a\in\M_+.$
\end{proposition}

\begin{proof}
Without loss of generality, we assume that $a\neq 0.$

First consider $\norm{a}\leq 1.$ By Proposition \ref{P01}, $(1-a)\in\M_+,$ so
\[\theta (a)\leq \theta (1-a)+\theta (a) =\theta (1-a+a)=\theta (1).\]

Now the conclusion follows if we apply this to $\dis \frac{a}{\norm{a}}$ instead of $a.$
\end{proof}

\begin{corollary}\label{C01} The following are equivalent for a weight $\theta$ on $\M_+$:
\begin{enumerate}
\item[(i)] $\theta(1)<\infty;$
\item[(ii)] $\theta(a)<\infty$ for all $a\in\M_+.$
\end{enumerate}\hfill $\Box$
\end{corollary}

A weight $\theta$ with $\theta (1)<\infty\;$ extends uniquely to $\widetilde{\theta} :\M\to\C$ by \[\widetilde{\theta}(a)=\big(\theta (b_1)-\theta (b_2)\big)+i \big(\theta (c_1)-\theta (c_2)\big),\] where $a$ has the unique decomposition mentioned in \S\ref{sec1}. Note that this canonical extension is a self-adjoint linear mapping, i.e. $\widetilde{\theta} (a^*) = \overline{\widetilde{\theta}(a)}$ for all $a\in\M.$

For convenience, we write $\theta$ instead of $\widetilde{\theta}$ whenever this unique extension is considered.

We recall that a weight is \begin{itemize}
             \item \emph{faithful}:  if $\theta(a)=0$ implies $a=0;$
             \item \emph{a state}: if $\theta(1)=1;$
             \item \emph{a trace}: if $\theta(a a^*)=\theta(a^* a)$ for all $a\in\M_;$ (equivalently, $\theta(u^* a u) =\theta (a)$ for all $a\in\M_+$ and all unitary $u\in\M,$ by \cite{Dix} I. Ch. 6 Cor.1);
             \item \emph{normal}: (\cite{Haa}) for any uniformly norm-bounded directed $\F\subset \M_+,$  if we have $\sup\F\in\M_+,\,$ then $\dis \theta\Big(\sup\F\Big)=\sup_{a\in\F}\theta(a);$
             \item $\kappa$-\emph{normal}: if the above holds for any uniformly norm-bounded directed $\F\subset \M_+$ with $\abs{\F} <\kappa.$
           \end{itemize}

For the rest of the section, $\M$ is assumed to be an internal $C^{*}$-algebra. Clearly, all the notions above have an internal counterpart. Moreover,  all the results obtained in the standard setting \textit{transfer} to the internal one.

We say that an internal weight $\theta :\M_+\to \s[0,\infty ]$ is \emph{$S$-continuous} if $\theta (\epsilon)\approx 0$ for all $0\approx\epsilon\in\M_+.$

\begin{lemma}\label{L01} Let $\theta :\M_+\to \s[0,\infty ]$ be an internal weight. Then the following are equivalent:
\begin{enumerate}
  \item [(i)] $\theta(1) <\infty.$
  \item [(ii)] $\theta$ is $S$-continuous.
  \item [(iii)] $\theta (a)<\infty$ whenever $a\in{\rm Fin}(\M_+).$
\end{enumerate}
\end{lemma}

\begin{proof}(i)$\Rightarrow$ (ii): Let $0\approx\epsilon\in\M_+.$ By Transfer of Proposition \ref{P02}, we have $\theta(\epsilon)\leq \theta(1)\norm{\epsilon}\approx 0,$ since $\theta(1) <\infty.$

(ii)$\Rightarrow$ (iii): Given $a\in{\rm Fin}(\M_+),$ suppose $\theta (a)\approx\infty.$ Then $\theta$ would have given the infinitesimal $\dis \frac{a}{\theta (a)}$ value $1,$ contradicting to (ii).

(iii)$\Rightarrow$ (i): Trivial.
\end{proof}

By a previous remark, we can then speak of an $S$-continuous weight $\theta$ on $\M,$ i.e. $\theta : \M\to \s\C.$

\begin{lemma}\label{L02}
Let $\theta$ be an $S$-continuous weight on $\M$ and $a, b\in {\rm Fin}(\M).$ Suppose $a\approx b,$ then $\theta (a)\approx \theta (b).$
\end{lemma}

\begin{proof} It suffices to prove that if $a\approx 0$ then $\theta(a)\approx 0.$ Since $a\approx 0$ can be decomposed as $a=a_1+i\,a_2,$ with $a_1\approx 0\approx a_2$ and $a_1, a_2\in \re{\M},$ we just need to prove that $\theta(a)\approx 0$ whenever $0\approx a\in\re{\M}.$ By \cite{Con} Chapter VIII, Proposition 3.4, $a$ can be decomposed as $a=a^\prime -a^{\prime\prime},$ where $a^\prime, a^{\prime\prime}\in\M_+$ are given by
\[a^{\prime} =f (a ),\; a^{\prime\prime} =g (a ),\quad \text{with }\, f(t)=\max(t,0)\; \text{and}\; g(t)=-\min(t,0).\]
(Here and in the sequel of the proof we follow the notation of \cite{Con}.)
Note that $f,g\in \CC \big(\sigma (a )\big)$ and  $\norm{f}_\infty, \norm{g}_\infty\approx 0,$ by $\sigma (a )\subseteq [-\norm{a},\norm{a}].$

By \cite{Con} Chapter VIII, Theorem 2.6 (b), we have $\norm{a^\prime}=\norm{f(a )}=\norm{f}_\infty$ and $ \norm{a^{\prime\prime}}=\norm{g(\epsilon )}=\norm{g}_\infty,$
therefore $a^\prime, a^{\prime\prime}\approx 0,$ and $\theta(a)\approx 0$ follows from $S$-continuity of $\theta.$
\end{proof}

As a consequence of Proposition \ref{P03} and Lemma \ref{L02}, we can define from an $S$-continuous internal weight $\theta$ on $\M$ a weight $\w\theta : \w\M \to \C$ by
$\w a\, \mapsto\,\st\theta(a).$

We want to prove that every  $S$-continuous internal weight in a $\kappa$-saturated nonstandard universe is $\kappa$-normal.  We begin with a definition and a remark.

Let $a,b\in\M.$ For $n\in\N^+,$ we define $a\leq_n b$ iff  there exist $a_1, b_1\in \M_+$ such that $a+a_1\leq b+b_1$ and $\norm{a_1},\norm{b_1}\leq n^{-1}.$

\begin{remark}\label{lleq} Let $\w a, \w b\in{\rm Fin}(\M).$ Then
 $$\w a\leq\w b\quad\Leftrightarrow\quad  a\leq_n b\ \text{\ for all \ }n\in\N^+.$$
The right-to-left implication is an application of Overspill. For the converse,
if $\w a\leq \w b,$ then $\epsilon=\sup\Big(\{\langle(a-b)h,h\rangle\, \vert\  h\in\HH_1\}\cup\{0\}\Big)$ is a nonnegative infinitesimal and $a\leq b+\epsilon\cdot 1.$ \hfill $\Box$
\end{remark}

\begin{theorem}\label{T01} Let $\theta :\M\to \s\C$ be an $S$-continuous weight in a $\kappa$-saturated nonstandard universe. Then $\w\theta$ is $\kappa$-normal.

\end{theorem}

\begin{proof}
Clearly $\w\theta$ is additive and positive homogeneous.

Let ${\mathcal A}\subseteq\w\M_+ $ be an  infinite norm-bounded directed family with $\abs{{\mathcal A}}<\kappa.$ Let $L$ be a norm-bound for the elements of ${\mathcal A}.$ Let  ${\mathcal A}_0$ be formed  by picking exactly one representative for each element in ${\mathcal A},$ so that ${\mathcal A}=\{\w a\,\vert a\in{\mathcal A}_0\}.$

Let $\, R= \sup\{\theta(\w a)\,\vert\  \w a\in{\mathcal A}\}.$ Since $\mathcal A$ is norm-bounded and $\theta$ is $S$-continuous, $R$ is finite (Proposition \ref{P02} and Lemma \ref{L01}).

We first claim that there exists $b\in {\rm Fin}(\M)$ such that
\[a\le_n b\quad\text{for all\  } a\in{\mathcal A}_0\text{\ and all\ }n\in\N^+\quad \text{and}\quad \theta(b)\approx R.\]

To see this, we let \[\F_{n,{\mathcal B}} = \Big\{ x\in\M\,\Big\vert\  a\leq_n x\ \text{\ for all\ } a\in{\mathcal B};\ \abs{\theta(x)-R}\leq\frac1n\text{\ and\ }\norm{x}\leq L+1\Big\},\]
where ${\mathcal B}\subset {\mathcal A}_0$ is finite and $n\in\N^+.$
By directedness of $\mathcal A$ and by Remark~\ref{lleq}, the $\F_{n,{\mathcal B}}$ are nonempty. Also, they have the finite intersection property, since $ \F_{\max(n,m), {\mathcal B}\cup {\mathcal C}}\subseteq \F_{n,{\mathcal B}}\cap \F_{m,{\mathcal C}}.$
By $\kappa$-saturation, we let $b\in \bigcap \F_{n,J},$ where the intersection ranges over  finite ${\mathcal B}\subset {\mathcal A}$ and $n\in\N^+.$ Then $b$ satisfies the claim above.

>From Remark~\ref{lleq}, we get $\w a\le\w b$ and hence $\w\theta(\w a)\le\w\theta(\w b)$ for all $\w a\in{\mathcal A}.$  Being $\w b$ a $\leq\,$-upper bound of ${\mathcal A},$ if $\sup{\mathcal A}$ exists then $$\sup\{\theta(\w a)\,\vert\  \w a\in{\mathcal A}\}\le\theta(\sup{\mathcal A})\le\w\theta(\w b)=\sup\{\theta(\w a)\,\vert\  \w a\in{\mathcal A}\}.$$
Therefore $\w\theta$ is $\kappa$-normal.
\end{proof}

\vskip 20pt

\section{Sequences and approximate properties}\label{app}

In this final section, we consider a standard $C^{*}$-algebra $\M$ and state some consequences of the previous results when the internal $C^{*}$-algebra is taken to be $\s\M$ - a nonstandard extension of $\M.$  Moreover, these consequences are stated without any references to nonstandard analysis nor to the nonstandard hull construction. They are examples of standard statements proved by means of the  nonstandard machinery.

We start with a simple observation: let $\mathcal N$ be an internal $C^*$-algebra and let $a\in\mathcal N$ satisfying $a\approx a^*\approx a^2.$ It follows  from the $C^*$-algebra identity  that $\norm{a}\approx \norm{a}^2,$ hence $a$ has finite norm.  Since $\w a=(\w a)^*=(\w a)^2$, $\w a$ is a projection and, by Theorem~\ref{plift01}(ii), there is some $p\in\proj{\mathcal N}$ such that $a\approx p$.

Now some definitions and notation. We call $\{ a_n\}_{n\in \N}\subset\M\,$ (not necessarily a convergent sequence) an \emph{almost projection sequence} (a.p.s.) if
\[\norm{a_n-a_n^*}+\norm{a_n^2-a_n}\,\to\, 0.\]
The sequence $\{ u_n\}_{n\in \N}\subset\M\,$ an \emph{almost partial isometry sequence} (a.p.i.s) if both $\{ u_n u_n^*\}_{n\in \N}\,$ and $\{ u_n^* u_n\}_{n\in \N}\,$ are a.p.s.

Let $\{ a_n\}_{n\in \N},\, \{ b_n\}_{n\in \N}\subset\M\,$ be a.p.s. Then we write $\{ a_n\}\sim_{\rm A} \{ b_n \}$ if there is an a.p.i.s. $\{ u_n\}_{n\in \N}\subset\M\,$ such that
\[\norm{u_n u_n^*-a_n}+\norm{u_n^* u_n -b_n}\,\to\, 0.\]

By the Transfer Principle, $\{ a_n\}_{n\in \N}\subset\M\,$ is an a.p.s. if and only if the internal extension $\{ a_n\}_{n\in \s\N}\subset\s\M\,$ satisfies the property that $a_n\approx a_n^*\approx a_n^2,$ for all infinite $n.$
Hence, by the above observation, we get that $\{ a_n\}_{n\in \N}\subset\M\,$ is an a.p.s. if and only if $\w a_n\in\proj {\widehat{\s\M}}$ for all infinite $n.$

Likewise, $\{ u_n\}_{n\in \N}\subset\M\,$ is an a.p.i.s. if and only if the internal extension $\{ u_n\}_{n\in \s\N}\subset\s\M\,$ satisfies the property that $\w u_n$ is a partial isometry in $\widehat{\s\M}\,$ for all infinite $n.$

Furthermore, for a.p.s. $\{ a_n\}_{n\in \N},\, \{ b_n\}_{n\in \N}\subset\M,\,$ we have $\{ a_n\}\sim_{\rm A} \{ b_n \}\,$ if and only if  $\w a_n \sim \w b_n\,$ for all infinite $n,$  where, as in Section~\ref{sec3},  $\sim$ denotes the Murray-von Neumann equivalence relation.

Using these characterizations, we apply the Transfer Principle to Theorem~\ref{plift01} and deduce the following:

\bigskip

\begin{corollary}\label{stanplift01} Let $\M$ be a standard $C^{*}$-algebra.

\begin{enumerate}
  \item [(i)]
  Let $\{a_n\}_{n\in \N}\subset\M\,$ be an a.p.s. Then for any $\epsilon\in\R^+\,$ there is $n\in\N\,$ such that
            \[\forall m\in\N\;\big[m>n\,\Rightarrow\,\exists p\in\proj{\M}\;\norm{a_m - p}\,<\,\epsilon\big].\]
 \item[(ii)] Let $\{ u_n\}_{n\in \N}\subset\M\,$ be an a.p.i.s. Then for any $\epsilon\in\R^+\,$ there is $n\in\N\,$ such that
      \[\forall m\in\N\;\big[m>n\,\Rightarrow\,\exists\;\text{partial isometry }u\in\M\, \norm{u_m - u}\,<\,\epsilon\big].\]
   \item [(iii)] Consider a.p.s. in $\M$ such that $\{ a_n\}\sim_{\rm A} \{ b_n\}.\,$ Then for any $\epsilon\in\R^+\,$ there is $n\in\N\,$ such that
      \[\forall m\in\N\;\big[ m>n\,\Rightarrow\,\exists p,q\in\proj{\M}\;\big(\,p\sim q\,\mbox{\ and\ }\,\norm{a_m-p}+\norm{b_m-q}\,<\,\epsilon\big)\big].\]
\end{enumerate}
\hfill $\Box$
\end{corollary}

\bigskip

From Theorem~\ref{finite} we have:

\bigskip

\begin{corollary}\label{stanfinite} Let $\M$ be a standard $C^{*}$-algebra.

\begin{enumerate}
  \item [(i)] If there is an a.p.s. $\{a_n\}_{n\in \N}\subset\M\,$ with $1\neq a_n$ for all $n\in\N\,$ and $\{ a_n\}\sim_{\rm A} \{ 1\},\,$ then $\M$ is infinite.
  \item[(ii)] If there are a.p.s. $\{ a_n\}_{n\in \N},\, \{ b_n\}_{n\in \N}\subset\M\,$ with
  \[a_nb_n\to 0\quad\text{and }\quad\{ a_n\}\sim_{\rm A} \{ 1\}\sim_{\rm A} \{ b_n\},\]
  then $\M$ is properly infinite.
\end{enumerate}
\hfill $\Box$
\end{corollary}

\bigskip

Note that the above corollaries can be stated for arbitrary infinite nets instead of countably infinite sequences. Moreover their converses are trivially true.

\bibliographystyle{amsplain}

\vfill

\small{\sc{\leftline{Dipartimento di Matematica} \leftline{Universit\`a
di Trento} \leftline{I-38050 Povo (TN)\quad Italy}

\noindent{\em baratell@science.unitn.it}

\vspace{20pt}

\leftline{School of Mathematical Sciences} \leftline{University of
KwaZulu-Natal} \leftline{Pietermaritzburg} \leftline{3209 South
Africa}

\noindent{\em ngs@ukzn.ac.za}}}

\vfill

\end{document}